\documentclass[%
bmf,%
 amsmath,amssymb,
preprint,%
]{revtex4-1}

\usepackage{graphicx}
\usepackage{dcolumn}
\usepackage{bm}

\begin{document}


\title[Nonresonant Hopf-Hopf bifurcation and a chaotic attractor in neutral functional differential equations]{Nonresonant Hopf-Hopf bifurcation and a chaotic attractor in neutral functional differential equations}

\author{Ben Niu}

\affiliation{Department of  Mathematics, Harbin Institute of Technology,\\  Harbin 150001, P.R. China.}%
\affiliation{Department of Applied Mathematics, Harbin University of Science and  Technology,\\ Harbin 150080, P.R. China.}%

\author{Weihua Jiang }

\affiliation{Department of  Mathematics, Harbin Institute of Technology,\\  Harbin 150001, P.R. China.}%


\begin{abstract}
 Nonresonant Hopf-Hopf singularity in neutral functional differential
equation (NFDE) is considered. An algorithm for calculating the
third-order normal form is established by using the formal adjoint
theory, center manifold theorem and the traditional normal form
method for RFDE. Van der Pol's equation with extended delay feedback
is studied as an example. The unfoldings near the Hopf-Hopf
bifurcation point is given by applying this algorithm. Periodic
solutions, quasi-periodic solutions  are found via theoretical
bifurcation diagram and numerical illustrations. The Hopf-Hopf
bifurcation diagram indicates the possible existence of a chaotic
attractor, which is confirmed by a sequence of simulations. [This is the full version of an article published as J. Math. Anal. Appl. 398 (2013) 362--371, doi:10.1016/j.jmaa.2012.08.051.]
\end{abstract}

\maketitle


\section{Introduction}
Normal form method is an important approach to the bifurcation
analysis, see \cite{Guckenheimer, Haleode, Wiggins2}. The key steps
for calculating the normal form of an ordinary differential equation
(ODE) are projecting original system onto the center manifold and
obtaining an approximate expression of normal form, up to any
desired degree of accuracy. In the case of functional differential
equations (FDEs), normal form is also an efficient method, but the
calculation process is quite tedious. Faria
 and Magalhaes's framework \cite{Faria1} obtained the normal forms for FDEs
 by recursive
changes of variables without computing beforehand the center
manifold of the singularity. Meanwhile, in the case of retarded
functional differential equations with parameters (RFDEs), where the
parameters are considered as new variables, the equation becomes an
abstract form of ODE in an enlarged phase space. Then by applying
the main idea in \cite{Faria1}, the computation of normal forms for
RFDEs with parameters are obtained in \cite{Faria2}. Another
approach for  calculating the normal form is due to Hassard {\it et
al} \cite{Hassard}, which mainly focused on the situation that Hopf
bifurcation appears. Both the two kinds of methods are widely used
in the bifurcation analysis of FDE. Normal form
method for  NFDE was developed recently, by \cite{Wee,Weedermann} in
which the author gave the computation procedure by employing the
method introduced by\cite{Faria1,Faria2}. In \cite{WangWei1}, normal
forms for NFDEs with parameters is established which is applied to
study the Hopf bifurcation in the  lossless transmission line, which
is the most famous equation of neutral type and has been studied a
lot, see \cite{WuBook, Halefde, WeiRuan, WuXia, KWu, ZW} and the
references cited therein.

Recent years, some authors turn to study the more complicated case,
i.e., codimension two bifurcations in RFDE.  Jiang {\it et al}
\cite{JiangYuan, JiangWang, WangJiang} studied the codimension two
bifurcations in a van der Pol's equation with nonlinear delay
feedback, which includes Bogdanov-Takens bifurcation,
Hopf-transcritical bifurcation, and Hopf-pitchfork bifurcation,
respectively. Zhang {\it et al} \cite{ZhangLi} studied the
Bogdanov-Takens bifurcation in a delayed predator-prey diffusion
system with a functional response. Ma, {\it et al} \cite{MaLu}
studied van der Pol's equation with a difference-type feedback,
where the feedback strength depends on the time lag. They mainly
studied the Hopf-Hopf bifurcation and observed some  interesting
phenomena such as stable tori, etc. However, most results about
codimension two bifurcations focus on the RFDE. Buono and
B$\acute{\textrm{e}}$lair \cite{buono} employed the methods
developed by \cite{Faria1,Faria2} to investigate the normal form and
universal unfolding of a vector field at non-resonant double Hopf
bifurcation points for particular classes of RFDEs. The case in NFDE
has not been studied.

In this paper, we extend the idea in \cite{Faria1, Faria2, Wee,
WangWei1} to the nonresonant Hopf-Hopf singularity in NFDE with
parameters:
\begin{equation}\label{NFDE}
    \frac d {dt}\left[Dx_t-G(x_t)\right]=L(\alpha)x_t+F(\alpha, x_t)
\end{equation}
where  $x_t\in C:=C([-\tau,0],R^n)$, $x_t(\theta):=x(t+\theta)$. $D$
and $L(\alpha)$ are bounded linear operators from $C$ to $R^n$ for
any $\alpha\in R^p$, with
$$D\phi=\phi(0)-\int_{-\tau}^0d[\mu(\theta)]\phi(\theta)$$ and
$$L(\alpha)\phi=\int_{-\tau}^0d[\eta(\theta,\alpha)]\phi(\theta)$$ for
$\phi\in C$, where $\mu(\theta)$ and $\eta(\theta,\alpha)$ are
matrix-valued functions of bounded variation which are continuous
from the left on $(-\tau,0)$ and such that $\eta(0,\alpha)=\mu(0)=0$
and $\mu$ is non-atomic at zero. Note that when $D\phi=\phi(0)$ and
$G(\phi)\equiv 0$, Eq.(\ref{NFDE}) degenerates to a RFDE. Thus the
method we established below is an extension to the RFDE case. Recall
that at a nonresonant (or \lq\lq no low-order resonant") Hopf-Hopf
bifurcation point, the corresponding characteristic equation has two
pairs of pure imaginary roots $\pm i\omega_+$ and $\pm i\omega_-$,
and further we have $\omega_-:\omega_+\neq 1:2~\textrm{or}~1:3$, if
we assume $0<\omega_-<\omega_+$.

The bifurcation results in NFDE are almost about codimension one
bifurcation such as Hopf bifurcation. As is known to all, studying
the codimension two bifurcation is a useful method to detect the
existence of homoclinic orbits, the coexistence of several periodic
orbits and the existence of quasi-periodic orbits (torus). More
precisely,  the universal unfoldings of the normal form is quite
important to reveal the dynamical behavior near the bifurcation
points. However, to our best knowledge, the universal unfoldings of
codimension two bifurcations in NFDE hasn't been well studied. The
current paper first employ the normal form method in FDE with
parameters given by \cite{Faria1, Faria2} to NFDE with Hopf-Hopf
singularity. For the sake of usage, we give an
 explicit and clear algorithm to deal with the Hopf-Hopf
singularity in detail.Firstly, the center manifold reduction and the
normal form derivation in the parameterized NFDE are presented.
After calculating the normal form near the Hopf-Hopf bifurcation
point, we show the dynamics of  NFDE near the critical point of the
Hopf-Hopf bifurcation is governed by a 3-dimensional system up to
the third order with unfolding parameters restricted on the center
manifold. Finally, it can be further reduced to a 2-dimensional
amplitude system, where these unfolding parameters can be expressed
by those perturbation parameters in the original NFDE. Our algorithm
is a  formulated procedure to study the dynamical behavior near a
nonresonant Hopf-Hopf bifurcation.

As an example to use these methods, we study the nonresonant
Hopf-Hopf bifurcation in van der Pol's equation with extended delay
feedback(See Pyragas\cite{Pyragas}), which is equivalent to a system
of NFDEs. Van der Pol's equation is widely studied by many authors
since it was first formulated for an electrical circuit with a
triode valve, for example \cite{Guckenheimer, Haleode,  Atay,
Maccari, MaLu, WeiJiang}. By analyzing the corresponding normal
form, we obtain the universal unfoldings near the Hopf-Hopf point.
Detailed bifurcation sets indicate the existence of  stable periodic
solution  and the stable quasi-periodic solution on torus. We find
that when the stable three-dimensional torus disappears, the van der
Pol's equation admits  a chaotic attractor.

The paper is organized as follows: in Section 2, we briefly present
the computation of normal forms for system (\ref{NFDE}) and give the
normal form derivation near the nonresonant Hopf-Hopf singularity.
Section 3 focuses on the van der Pol's equation. The conditions of
the existence of Hopf-Hopf bifurcation is obtained. Using the
algorithm presented in Section 2, the corresponding normal form is
calculated, and the detailed bifurcation sets are drawn. Appropriate
simulations are carried out to illustrated the theoretical results.
Finally a conclusion section is given, in which we also give some
discussions.

\section{Reduction and normal form for NFDEs with Hopf-Hopf Singularity}
In this section, we present the regular normal form method for
system (\ref{NFDE}), and then calculate the normal form near a
nonresonant  Hopf-Hopf bifurcation point.
\subsection {Normal form derivation in NFDE}
In this paper we always
 assume that $x_t$ is differentiable, $F$ and $G$
are $C^N$-smooth, $N\geq 3$, $F(\alpha, 0)=G(0)=0$, $F'(\alpha,
0)=G'(0)=0$ and $G$ doesn't depend on $\phi(0)$ in (\ref{NFDE}),
where the notation $'$ stands for the Fr$\acute{\textrm{e}}$chet
derivative. Under these hypothesis, obviously $x_t=0$ is an
equilibrium of (\ref{NFDE}) which is equivalent to
\begin{equation}\label{NFDEa}
    \frac d {dt}Dx_t=L(\alpha)x_t+F(\alpha, x_t)+G'(x_t)\dot x_t
\end{equation}
By introducing the enlarged phase space $BC$ in which the functions
from $[-\tau,0]$ to $R^n$ are uniformly continuous on $[-\tau,0)$
with a possibly jump discontinuously at $0$, Eq.(\ref{NFDEa}) can be
written as an abstract ordinary differential equation on $BC$:
\begin{equation}\label{NFDEabstract}
\frac d {dt}
x_t=Ax_t+X_0[L(\alpha)-L(0)]x_t+X_0F(\alpha,x_t)+X_0G'(x_t)\dot x_t
\end{equation}
where
\begin{equation}\label{A}A\phi:=\phi'+X_0[L(0)\phi-D\phi'].\end{equation}
$A$ is the infinitesimal generator of the semigroup of solutions to
the linear system $$\frac d {dt}Dx_t=L(0)x_t.$$ $X_0(\theta)=0$ for
$-\tau\leq \theta<0$ and $X_0(0)=Id_{n\times n}$. Introducing
$\alpha$ as a new variable, we have
\begin{equation}\label{NFDEnonpara}
\begin{array}{rcl}\frac d {dt}
x_t&=&Ax_t+X_0[L(\alpha)-L(0)]x_t\\&&+X_0F(\alpha,x_t)+X_0G'(x_t)\dot
x_t\\ \frac{d}{dt}\alpha(t)&=&0\end{array}
\end{equation}
which can be considered as an ODE  with no parameters in the product
space $\widetilde{{BC}}:=BC\times R^p$. Following \cite{Faria1,
Faria2, Wee} we summarize the calculation of the normal forms for
Eq.(\ref{NFDEnonpara}) as follows. Noting that here we use $x_t\in
C:=C([-\tau,0],R^n)$ as mentioned in section 1. If we use complex
vectors to decompose the phase space, these discussions also hold
true for the complex case $x_t\in C:=C([-\tau,0],\mathbb{C}^n)$ when
the operators $L$, $D$, $F$ are extended to complex functions in the
natural way.

Decompose $\widetilde{{BC}}$ by
$\widetilde{{BC}}=\widetilde{P}\bigoplus
\textrm{Ker}\widetilde{\pi}$, where $\widetilde{P}=P\times R^p$, $P$
is the generalized eigenspace for $A$ associated with a nonempty
finite set $\Lambda$ of eigenvalues of $A$. $\widetilde{\pi}$ is the
projection of $\widetilde{BC}$ upon $\widetilde{P}$.
$\Phi=(\phi_1,\phi_2,\cdots,\phi_m)$ is a basis for $P$ with
$(\Psi,\Phi)=Id_{m\times m}$ where
$\Psi=(\psi_1,\psi_2,\cdots,\psi_m)$ is a basis for $P^\ast$, the
dual space of $P$.  The bilinear form $(\cdot,\cdot)$ is
\begin{equation}\label{bilinear}\begin{array}{l}
    (\psi,\phi)=\psi(0)\phi(0)-\int_{-\tau}^0d\left[\int_0^\theta
    \psi(\xi-\theta)d\mu(\xi)\right]\phi(\theta)\\+\int_{-\tau}^0\int_0^\theta\psi(\xi-\theta)d\eta(\theta,0)\phi(\xi)d\xi\end{array}
\end{equation}
Choose $B$ such that $A\Phi=\Phi B$. If we decompose $$\left(
              \begin{array}{c}
                x_t \\
                \alpha_t \\
              \end{array}
            \right)=\left(
                      \begin{array}{cc}
                        \Phi & 0 \\
                        0 & Id_{p\times p} \\
                      \end{array}
                    \right)
            \left(
                      \begin{array}{c}
                        z(t) \\
                        \alpha(t) \\
                      \end{array}
                    \right)+\left(
                              \begin{array}{c}
                                w_1 \\
                                w_2 \\
                              \end{array}
                            \right),
$$
where $(z(t),\alpha(t))\in R^{m+p}$ and $(w_1,w_2)\in
\textrm{Ker}(\widetilde{\pi}) $, then (\ref{NFDEnonpara}) is
equivalent to
\begin{equation}\label{decomposition}
    \begin{array}{l}\left(
       \begin{array}{c}
         \dot z \\
         \dot \alpha \\
       \end{array}
     \right)=\left(
               \begin{array}{c}
                 Bz \\
                 0\\
               \end{array}
             \right)+\left(
                       \begin{array}{c}
                         \Psi(0)(L(\alpha(0)+w_2(0))-L(0))(\Phi z+w_1) \\
                         0 \\
                       \end{array}
                     \right)\\+\left(
                               \begin{array}{c}
                                \Psi(0)F(\Phi z+w_1, \alpha(0)+w_2(0))\\
                                 0 \\
                               \end{array}
                             \right)\\+\left(
                                       \begin{array}{c}
                                         \Psi(0)G'(\Phi z+w_1)(\Phi\dot z+\dot w_1)\\
                                         0 \\
                                       \end{array}
                                     \right)\\
                                     \left(
       \begin{array}{c}
         \dot w_1 \\
         \dot w_2 \\
       \end{array}
     \right)=\left(
               \begin{array}{c}
                 A_{Q^1}w_1 \\
                 \dot w_2-Y_0\dot w_2(0)\\
               \end{array}
             \right)\\+\left(
                       \begin{array}{c}
                         (Id-\pi)X_0(L(\alpha(0)+w_2(0))-L(0))(\Phi z+w_1) \\
                         0 \\
                       \end{array}
                     \right)\\+\left(
                               \begin{array}{c}
                                (Id-\pi)X_0F(\Phi z+w_1, \alpha(0)+w_2(0))\\
                                 0 \\
                               \end{array}
                             \right)\\+\left(
                                       \begin{array}{c}
                                         (Id-\pi)G'(\Phi z+w_1)(\Phi\dot z+\dot w_1)\\
                                         0 \\
                                       \end{array}
                                     \right)\end{array}
\end{equation}
where the newly defined $A_{Q^1}$ is the restriction of $A$ to
$Q^1:=Q\bigcap C^1$ with $Q$ being the complementary space of $P$ in
$C$. $Y_0(\theta)=0$ for $-\tau\leq \theta<0$ and
$Y_0(0)=Id_{m\times m}$. Noting that $w_2(0)=0$ because $w_2\in
R^p$, and dropping the auxiliary equations we get the equation
(\ref{decomposition}) in $BC=P\bigoplus \textrm{Ker}(\pi)$
equivalent to
\begin{equation}\label{decomsimple}\begin{array}{l}
    \dot z=Bz + \Psi(0)[(L(\alpha)-L(0))(\Phi z+w_1)\\~~~~~~~~+F(\Phi
    z+w_1,\alpha)+G'(\phi z+w_1)(\Phi \dot z+\dot w_1)]\\ \dot
    w_1=A_{Q^1}w_1+(Id-\pi)X_0[(L(\alpha)-L(0))(\Phi z+w_1)\\~~~~~~~~+F(\Phi
    z+w_1,\alpha)]\\~~~~~~~~+(Id-\pi)G'(\Phi z+w_1)(\Phi \dot z+\dot
    w_1)\end{array}
\end{equation}
Write the Taylor expansion of (\ref{decomsimple}), we have
\begin{equation}\label{taylor}\begin{array}{lll}
    \dot z&=&Bz+\sum_{j\geq 2}\frac 1{j!}f_j^1(z,w_1,\alpha)\\
    \dot w_1&=&A_{Q^1}w_1+\sum_{j\geq 2}\frac
    1{j!}f_j^2(z,w_1,\alpha)\end{array}
\end{equation}
To derive the normal form of the $j$th order, we make the
transformations of variables for $j\geq 2$, $(z,w_1,\alpha)\mapsto
(\hat{z},\hat{w}_1,\hat{\alpha})$ given by
\begin{equation}\label{transform}
(z,w_1,\alpha)=(\hat{z},\hat{w}_1,\hat{\alpha})+\frac
1{j!}\widetilde{U_j}(\hat{z}, \hat{\alpha})
\end{equation}
 with
$\widetilde{U_j}=({U_j^1},{U_j^2},{U_j^3})\in V_j^{m+p}(R^m)\times
V_j^{m+p}(Q^1)\times V_j^{m+p}(R^p)$, $U_j=(U_j^1,U_j^2)$, where for
a normed space $X$, we denote by $V_j^{m+p}(X)$ the linear space of
homogeneous polynomials of degree $j$ in $m+p$ real variables with
coefficients in $X$. To compute the normal form we define the
operator $M_j$ on $V_j^{m+p}(R^m\times Ker\pi)$ by
$M_j(q,h)=(M_j^1q,M_j^2h)$, where
$$(M_j^1q)(z,\alpha)=D_zq(z,\alpha)Bz-Bq(z,\alpha),$$
$$(M_j^2h)(z,\alpha)=D_zh(z,\alpha)Bz-A_{Q_{bt}^1}h(z,\alpha),$$
with $q(z,\alpha)\in V_j^{m+p}(R^m)$, $h(z,\alpha)(\theta)\in
V_j^{m+p}(Q^1)$. Then we have the following decompositions
$$V_j^{m+p}(R^m)=Im(M_j^1)\bigoplus Im(M_j^1)^c,$$
$$V_j^{m+p}(R^m)=Ker(M_j^1)\bigoplus Ker(M_j^1)^c,$$
$$V_j^{m+p}(Ker\pi)=Im(M_j^2)\bigoplus Im(M_j^2)^c,$$
$$V_j^{m+p}(Q^1)=Ker(M_j^2)\bigoplus Ker(M_j^2)^c.$$ Denote the
projections associated with the above decompositions of
$V_j^{m+p}(R^m)\times V_j^{m+p}(Ker\pi)$ over $Im(M_j^1)\times
Im(M_j^2)$ and of $V_j^{m+p}(R^m)\times V_j^{m+p}(Q^1)$ over
$Ker(M_j^1)^c\times Ker(M_j^2)^c$ by, respectively
$P_{I,j}=(P_{I,j}^1,P_{I,j}^2) $, and
$P_{K,j}=(P_{K,j}^1,P_{K,j}^2)$. By transformation
(\ref{transform}), the $j$th order term in the normal form becomes
$g_j=\bar f_j-M_jU_j$, where $\bar f_j$ denotes the terms of order
$j$ obtained after computation of the normal form up to order $j-1$.
Following \cite{Faria1, Faria2} we have an adequate choice of $U_j$
by
$$U_j(z,\alpha)=M_j^{-1}P_{I,j}\bar f_j(z,0,\alpha)$$ and thus $g_j(z,0,\alpha)=(I-P_{I,j})\bar
f_j(z,0,\alpha)$.

Following the general work in \cite{Faria1,Faria2}, and using the
center manifold theory presented in \cite{ Hassard, Wiggins1, Chow,
Carr} we have the following conclusion:

{\bf Theorem 1} Suppose that in system (\ref{NFDE}) the
infinitesimal generator $A$ has $m$ eigenvalues with zero real
parts, and the other eigenvalues have negative real parts. Denote
$\Lambda=\{\lambda | \lambda\in\sigma(A)~ \textrm{and
~Re}\lambda=0\}$, and the corresponding generalized eigenspace
spanned by $\Phi=(\phi_1,\phi_2,\cdots,\phi_m)$ with $A\Phi=\Phi B$.
Assume further that the nonresonance conditions (Ref.
{\cite{Faria2}}) relative to $\Lambda$ are satisfied. Then the
dynamics in (\ref{NFDE}) near $x_t=0$ are governed by
\begin{equation}\label{NFALL}
    \dot z=Bz+\sum_{j\geq 2}\frac{1}{j!}g_j(z,0,\alpha)
\end{equation}

\subsection{Normal form of Hopf-Hopf bifurcation}

Generally,  Hopf-Hopf bifurcation occurs in Eq.(\ref{NFDE}) when
$\alpha=(\alpha_1,\alpha_2)=0$ if in $\sigma(A)$ there are four
points with zero real parts, $\{\pm i \omega_1, \pm i \omega_2\}$.
This is just to say the characteristic equation of Eq.(\ref{NFDE}),
\begin{equation}\label{CENFDE}
    det(\Delta(\lambda))=det(\lambda D(e^{\lambda
    \cdot})-L(e^{\lambda\cdot}))=0
\end{equation}
has two pairs of pure imaginary roots. Without loss of generality,
we assume $\omega_1<\omega_2$. We decompose $x_t=\Phi z(t)+w$, with
$z(t)=(z_1(t),z_2(t),z_3(t), z_4(t))^T\in \mathbb{C}^4$ and $w\in
\textrm{Ker}(\pi)$. Following the method in section 2.1, we first
calculate $\Phi(\theta),\theta\in[-\tau,0]$ and
$\Psi(s),s\in[0,\tau]$ which satisfy $A\Phi=\Phi B$, with
$$B=\left(
                                                 \begin{array}{cccc}
                                                  i\omega_1 & 0 & 0 &  0 \\
                                                  0 & -i\omega_1  & 0 & 0 \\
                                                  0 & 0 & i\omega_2 & 0 \\
                                                  0 & 0 & 0 & -i\omega_2 \\
                                                 \end{array}
                                               \right)
,$$ $A^\ast \Psi=B\Psi$, with
\begin{equation}\label{Aast}A^{\ast}\nu
=-\frac{d\nu}{ds},\end{equation} where $Dom(A^{\ast})=\{\nu\in
C^{\ast}:=C([0,\tau],\mathbb{C}^{n\ast}), \frac{d\nu}{ds}\in
C^{\ast}, D\frac{d\nu}{ds}=-\int_0^{\tau}\nu(s)d[\eta(-s,0)]\},$
and $(\Psi,\Phi)=Id_{3\times 3}.$ Here $\mathbb{C}^{n\ast}$ is the
$n$ dimensional space with row vectors. Noting that, compared with a
RFDE, the operator $D$ only changes the definition $A$ and $A^\ast$.

Recall that $\omega_1:\omega_2\neq 1:2, ~\textrm{or}~1:3$ we have
$(Im(M_2^1))^c$ is spanned by the elements
$$\{z_1\alpha_ie_1,~z_2\alpha_ie_2,~z_3\alpha_ie_3,~z_4\alpha_ie_4\},~i=1,2,$$
with
$e_1=(1,0,0,0)^T,~e_2=(0,1,0,0)^T,~e_3=(0,0,1,0)^T,~e_4=(0,0,0,1)^T$.
Thus the normal form of (\ref{NFDE}) on the center manifold of the
origin near $(\alpha_1,\alpha_2)=0$ has the form
\begin{equation}\label{NF2}
    \dot z=Bz+\frac 1 2 g_2^1(z,0,\alpha)+h.o.t.,
\end{equation}
with $g_2^1(z,0,\alpha)=Proj_{(Im(M_2^1))^c}f_2^1(z,0,\alpha).$

To find the third order normal form of the Hopf-Hopf singularity,
let $M_3$ denote the operator defined in $V_3^3(\mathbb{C}^3\times
\textrm{Ker}(\pi))$. Here we neglect the high order term of the
perturbation parameters $(\alpha_1,\alpha_2)$. Thus $(Im(M_3^1))^c$
is spanned by
$$\{z_1^2z_2e_1, z_2^2z_1e_2, z_3^2z_4e_3, z_4^2z_3e_4, z_1z_3z_4e_1, z_2z_3z_4e_2, z_1z_2z_3e_3, z_1z_3z_4e_4\}.$$

The normal form of (\ref{NFDE}) up to the third order is
\begin{equation}\label{NF3}
\dot z=Bz+\frac 1 {2!} g_2^1(z,0,\alpha)+\frac 1 {3!}
g_3^1(z,0,0)+h.o.t.,
\end{equation}
where $g_3^1(z,0,0)=Proj_{(Im(M_3^1))^c}\bar f_3^1(z,0,0),$ with
$$(\bar f_3^1, \bar f_3^2)^T=( f_3^1, f_3^2)^T+\frac 32[(D_{z,
w}(f_2^1,f_2^2)^TU_2-(D_{z, w}U_2)(g_2^1,g_2^2))]$$ and
$$U_2(z,\alpha)=(U_2^1, U_2^2)^T=M_2^{-1}P_{I,2}f_2(z,0,\alpha).$$

Recall that system (\ref{NFDE}) undergoes a Hopf-Hopf bifurcation at
$x_t=0$ when $\alpha_1=\alpha_2=0$. Assume further all the other
roots except $\Lambda=\{\pm i\omega_1,\pm i\omega_2\}$ have negative
real parts, which obviously means that the nonresonance conditions
relative to $\Lambda$ are satisfied. Following Theorem 1, we have
that the dynamical behavior of (\ref{NFDE}) near $x_t=0$ is governed
by the general normal form of the third order
\begin{equation}\label{NF23}\begin{array}{lll}
    \dot
    z_1&=&i\omega_1z_1+a_{11}\alpha_1z_1+a_{12}\alpha_2z_1+c_{11}z_1^2z_2+c_{12}z_1z_3z_4,\\
    \dot
    z_2&=&-i\omega_1z_2+\bar a_{11}\alpha_1z_2+\bar a_{12}\alpha_2z_2+\bar c_{11}z_1z_2^2+\bar c_{12}z_2z_3z_4,\\
    \dot
    z_3&=&i\omega_2z_3+a_{21}\alpha_1z_3+a_{22}\alpha_2z_3+c_{21}z_1z_2z_3+c_{22}z_3^2z_4,\\
    \dot
    z_4&=&-i\omega_2z_4+\bar a_{21}\alpha_1z_4+\bar a_{22}\alpha_2z_4+\bar c_{21}z_1z_2z_4+\bar c_{22}z_3z_4^2,\end{array}
\end{equation}
Make the transformation
$z_1=r_1\cos\theta_1+ir_1\sin\theta_1,~z_2=r_1\cos\theta_1-ir_1\sin\theta_1,$
$z_3=r_2\cos\theta_2+ir_2\sin\theta_2,~z_4=r_2\cos\theta_2-ir_2\sin\theta_2,$,$r_1,r_2>0,$
then we have the amplitude system

\begin{equation}\label{NF23polar}\begin{array}{lll}
    \dot
   r_1&=&\textrm{Re}a_{11}\alpha_1r_1+\textrm{Re}a_{12}\alpha_2r_1+\textrm{Re}c_{11}r_1^3+\textrm{Re}c_{12}r_1r_2^2,\\
    \dot
   r_2&=&\textrm{Re}a_{21}\alpha_1r_2+\textrm{Re}a_{22}\alpha_2r_2+\textrm{Re}c_{21}r_1^2r_2+\textrm{Re}c_{22}r_2^3,\end{array}
\end{equation}

Denote by $\epsilon_1=\textrm{Sign}(\textrm{Re}c_{11})$,
$\epsilon_2=\textrm{Sign}(\textrm{Re}c_{22})$.  After re-scaling
$\hat r_1=r_1 \sqrt{|\textrm{Re}c_{11}|}$, $\hat r_2=r_2
\sqrt{|\textrm{Re}c_{22}|}$ and $\hat t=t \epsilon_1$, then
Eq.(\ref{NF23polar}) becomes, after dropping the hats

\begin{equation}\label{NF23polarreduce}\begin{array}{lll}
       \dot r_1&=&r_1(c_1+r_1^2+b_0r_2^2)+h.o.t.,\\ \dot r_2&=& r_2(c_2+c_0r_1^2+d_0r_2^2)+h.o.t.,
\end{array}
\end{equation}
where
\begin{eqnarray}\label{array}
 \nonumber c_1 &=& \epsilon_1\textrm{Re}a_{11}\alpha_1+\epsilon_1\textrm{Re}a_{12}\alpha_2 \\
 \nonumber c_2 &=& \epsilon_1\textrm{Re}a_{21}\alpha_1
+\epsilon_1\textrm{Re}a_{22}\alpha_2 \\
  b_0 &=& \frac{\epsilon_1\epsilon_2\textrm{Re}c_{12}}{\textrm{Re}c_{22}} \\
  \nonumber c_0 &=& \frac{\textrm{Re}c_{21}}{\textrm{Re}c_{11}} \\
  \nonumber d_0 &=& \epsilon_1\epsilon_2
\end{eqnarray}
Applying the results in section 7.5 of \cite{Guckenheimer}, Eq.
(\ref{NF23polarreduce}), truncated up to the third order,  has
twelve distinct types of unfoldings with respect to different signs
of $b_0,~c_0,~d_0$, and $d_0-b_0c_0$, which is shown in Table 1. The
detailed phase portraits can be found in  \cite{Guckenheimer}. Here
we only state the VIa case for the sake of usage.

When $b_0>0,~c_0<0,~d_0=-1$ and $d_0-b_0c_0>0$, case VIa arise. Near
the bifurcation point, the $\alpha_1$--$\alpha_2$ plane is divided
by eight lines:
\begin{description}
\item [$L_1$]: $c_2=0,~c_1>0$;
\item [$L_2$]: $c_1=0,~c_2>0$;
\item [$L_3$]: $c_2=c_0c_1,~c_2>0$;
\item [$L_4$]: $c_2=\frac{c_0-1}{b_0+1}c_1+O(c_1^2),~c_2>0$;
\item [$L_5$]: $c_2=\frac{c_0-1}{b_0+1}c_1,~c_2>0$;
\item [$L_6$]: $c_2=-\frac{c_1}{b_0},~c_2>0$;
\item [$L_7$]: $c_2=0,~c_1<0$;
\item [$L_8$]: $c_1=0,~c_2<0$;
\end{description}
Eight different phase portraits, when parameters lie between every
two neighboring lines, are list in Figure \ref{fig3} on page
\pageref{fig3}.

 {\centering {{\bf Table 1.} The twelve unfoldings of system
(\ref{NF23polarreduce}).}\\
\begin{tabular}{lllllllllllll}
  \hline
  Case & Ia & Ib & II & III & IVa & IVb & V & VIa & VIb & VIIa & VIIb & VIII \\
  \hline
  $d_0$ & +1 & +1 & +1 & +1 & +1 & +1 & --1 & --1 & --1 & --1 & --1 & --1 \\
   $b_0$ & + & + & + & -- & -- &-- & +& + & + & -- & -- & -- \\
   $c_0$ & + & + & -- & + & -- & -- & + & -- & -- & + & +& -- \\
   $d_0-b_0c_0$ & + & -- & + & + &+ & -- & -- & + & -- & + &-- & -- \\
  \hline
\end{tabular}}

So far, we have given the whole algorithm to determine the
unfoldings in a NFDE with Hopf-Hopf singularity, which includes
three key steps:
\begin{description}
  \item[\textit{Step 1}] Analyzing the associated characteristic equation to
  obtain the condition under which a Hopf-Hopf bifurcation
  occurs.
  \item[\textit{Step 2}] Write the equivalent ODE in the enlarged phase
  space. Calculate $\Phi$ and $\Psi$ by Eq.(\ref{A}) and (\ref{Aast}).
  \item[\textit{Step 3}] Decomposing the original system as
  Eq.(\ref{decomsimple}) and finally obtaining the normal form
as Eq.(\ref{NF23}). Calculate $b_0,~c_0,~d_0$, and $d_0-b_0c_0$.
\end{description}

\section{Hopf-Hopf bifurcation in van der Pol's equation with extended delay feedback}
In this section van der Pol's equation with extended delay feedback
is studied. Hopf-Hopf points are detected by analyzing the
associated characteristic equation. Near these points, we calculate
the normal form  by the  algorithm given in Section 2, and all the
key values are obtained. A numerical example  provides several kinds
of interesting phenomena which illustrated the theoretical results
given by bifurcation sets.
\subsection{The existence and the normal form derivation}
In this section we will study the Hopf-Hopf bifurcation in van der
Pol's equation with extended delay feedback, basing on the method
presented in section 2.

Consider the following van der Pol's equation
 \begin{equation} \label{vdpa}\ddot x+\varepsilon (x^2-1)\dot x+x=\varepsilon k
\vartheta(t)\end{equation} where $\varepsilon>0$. $k$ is the
strength of the feedback $\vartheta(t)$, which is the linear part in
the feedback signal. $\vartheta(t)$ depends on the current state and
a sequence of the past states, which is defined by
\begin{equation}\label{vdpb}\vartheta(t)=(1-\mu)x(t)+\mu
\vartheta(t-\tau),\end{equation} with $0<\mu<1$. Eq. (\ref{vdpa}) is
equivalent to
 \begin{equation}
\label{vdpc}\mu\ddot x(t-\tau)+\mu\varepsilon (x^2(t-\tau)-1)\dot
x(t-\tau)+\mu x(t-\tau)=\mu\varepsilon k
\vartheta(t-\tau)\end{equation} Using (\ref{vdpa})-(\ref{vdpc}), we
have
\begin{equation}\label{vdpd}\begin{array}{l}\ddot x-\mu\ddot x(t-\tau)+\varepsilon (x^2-1)\dot
x-\mu\varepsilon (x^2(t-\tau)-1)\dot x(t-\tau)+x-\mu
x(t-\tau)\\~~~~~~~~~~~~~~=\varepsilon
k(1-\mu)x(t)\end{array}\end{equation} This is a NFDE of second
order.
Introduce a new variable $y(t)=\dot x(t)$, then (\ref{vdpd}) becomes a system of NFDEs  \begin{equation}\label{vdp}\left\{\begin{array}{l}\dot x=y\\
\dot y-\mu\dot y(t-\tau)=[-1+\varepsilon k(1-\mu)]x+\varepsilon
y+\mu x(t-\tau)-\varepsilon\mu
y(t-\tau)\\~~~~~~~~~~~~~~~~~~~~~~~~~~~-\varepsilon
x^2y+\varepsilon\mu
x^2(t-\tau)y(t-\tau)\end{array}\right.\end{equation}

We begin with the trivial equilibrium  $E_0=(0,0)^T$ of (\ref{vdp}).
The  characteristic equation of the corresponding  linearized
equation is
\begin{equation}\label{CE}\lambda^2-\mu\lambda^2e^{-\lambda\tau}-\varepsilon\lambda+\varepsilon\mu\lambda
e^{-\lambda\tau}-\mu e^{-\lambda\tau}+1-\varepsilon
k(1-\mu)=0.\end{equation}

Now, we start analyzing the Hopf-Hopf bifurcation in
(\ref{vdp}) following the three steps state in the previous section.\\
\textbf{\textit{Step 1.}} We study the existence
 of the Hopf-Hopf bifurcation via detect the interjection of the Hopf bifurcation curves in Eq.(\ref{vdp}).

To detect the conditions that Hopf bifurcation occurs in
$(\ref{vdp})$, we substitute $\lambda=i\omega$, $\omega>0$ into
$(\ref{CE})$. Separating the real and imaginary parts gives
\begin{equation}\label{sepa}\left\{\begin{array}{l}(\mu\omega^2-\mu)\cos\omega\tau+\varepsilon\mu\omega\sin\omega\tau=\omega^2-1+\varepsilon
k(1-\mu)\\-(\mu\omega^2-\mu)\sin\omega\tau+\varepsilon\mu\omega\cos\omega\tau=\varepsilon\omega\end{array}\right.\end{equation}
which solves \begin{equation}\label{cossin}\left\{\begin{array}{l}
\cos(\omega\tau)=\frac{(\mu\omega^2-\mu)(\omega^2-1+\varepsilon
k(1-\mu))+(\varepsilon \mu\omega)(\varepsilon\omega)}
{(\mu\omega^2-\mu)^2+(\varepsilon\mu\omega)^2}\\
\sin(\omega\tau)=\frac{-(\mu\omega^2-\mu)(\varepsilon\omega)+(\varepsilon\mu\omega)(\omega^2-1+\varepsilon
k(1-\mu))}{(\mu\omega^2-\mu)^2+(\varepsilon\mu\omega)^2}\end{array}\right.\end{equation}
Hence, we
have\begin{equation}\label{ome2}(\mu\omega^2-\mu)^2+(\varepsilon\mu\omega)^2=(\omega^2-1+\varepsilon
k(1-\mu))^2+(\varepsilon\omega)^2\end{equation} which is equivalent
to
\begin{equation}\label{rho}W(\rho)=a\rho^2+b\rho+c=0\end{equation} where $\rho=\omega^2$,
$a=(1+\mu)$, $b=[2\varepsilon k-2(1+\mu)+\varepsilon^2(1+\mu)]$,
$c=\varepsilon^2k^2(1-\mu)-2\varepsilon k+1+\mu$.

Assume
$$(H1):k<\min\left\{\frac 1\varepsilon,
\frac{1+\mu}{\varepsilon}-\frac{\varepsilon(1+\mu)}{2}\right\},$$ we
have $c>0, b<0$ and Lemma 2.1 holds.  Further more if
$$(H2):\Delta=(b^2-4ac)>0,$$
then (\ref{ome2}) solves by two positive roots
$\omega_\pm=\sqrt{\rho_\pm},$ where
$\rho_\pm=\frac{-b\pm\sqrt{b^2-4ac}}{2a}$. Denote by $\tau_0^+$ (or
$\tau_0^-$) the unique root of Eq.(\ref{cossin}) when
$\omega=\omega_+$ (or $\omega=\omega_-$), such that
$\omega\tau_0^\pm \in [0,2\pi)$. Also denote by
\begin{equation}\label{tauj}\tau_j^\pm=\tau_0^\pm+\frac{2j\pi}{\omega_\pm},~j=0,1,2,\cdots\end{equation}

Take the derivative with respect to $\tau$ in Eq.(\ref{CE}), and use
Eq.(\ref{cossin}). After a few straightforward calculations, we have
\begin{equation}\label{transcon}
\textrm{Sign}\left(\left.\textrm{Re}\frac{d\lambda}{d\tau}\right|_{\tau=\tau_j^\pm}\right)=\textrm{Sign}(\left.W'(\rho))\right|_{\rho=\omega_\pm^2}.\end{equation}

Base on the above preparation, together with the Hopf bifurcation
theorem in \cite{WeiRuan}, we can give the conclusions about Hopf
bifurcation in $(\ref{vdp})$.

{\bf Theorem 2} Consider system $(\ref{vdp})$, with $\varepsilon>0$.
Assume $(H1)$, $(H2)$ hold.  If $\tau_0^->\tau_0^+$, then $E_0$ of
system (\ref{vdp}) is unstable for any $\tau\geq 0$. If
$\tau_0^-<\tau_0^+$ then there exists an integer $m\geq 0$ such that
$E_0$ is stable when $\tau\in(\tau_0^-, \tau_0^+)\cup(\tau_1^-,
\tau_1^+)\cup\cdots\cup(\tau_m^-,\tau_m^+)$, and is unstable when
$\tau\in(0, \tau_0^-)\cup(\tau_0^+,
\tau_1^-)\cup\cdots\cup(\tau_m^+,+\infty)$. Moreover, system
(\ref{vdp}) undergoes a  Hopf bifurcation at $\tau_j^+$(or
$\tau_j^-$), $j=0,1,2,\cdots$.

Now, we're in position to give an existence condition that a
Hopf-Hopf bifurcation occurs. Basing on the preparation about Hopf
bifurcation, we detect the possible existence of the Hopf-Hopf
bifurcation point, which is the interjection of two Hopf bifurcation
curve. If we fix $\varepsilon$ and $\mu$, then a figure like Figure
\ref{fig1} is drawn, in which the points denoted by $HH1$ and $HH2$
are Hopf-Hopf bifurcations. The exactly critical value can be obtain
by the calculation process:
\begin{description}
  \item[1] Solve $\omega_\pm$ as the function of $k$ from
  Eq.(\ref{rho}).
  \item[2] Substitute $\omega_\pm$ and $\tau_j^\pm$ into
  Eq.(\ref{cossin}). For $j=j_0$, solve $k=k_0$ from \begin{equation}\label{k0}\begin{array}{l}\left(\arccos\frac{(\mu\omega_+^2-\mu)(\omega_+^2-1+\varepsilon
k(1-\mu))+(\varepsilon \mu\omega_+)(\varepsilon\omega_+)}
{(\mu\omega_+^2-\mu)^2+(\varepsilon\mu\omega_+)^2}+2j\pi\right)/\omega_+\\~~~~~~~~~=\left(\arccos\frac{(\mu\omega_-^2-\mu)(\omega_-^2-1+\varepsilon
k(1-\mu))+(\varepsilon \mu\omega_-)(\varepsilon\omega_-)}
{(\mu\omega_-^2-\mu)^2+(\varepsilon\mu\omega_-)^2}+2j\pi\right)/\omega_-.\end{array}\end{equation}
  \item[3] Compute $\tau_0=\tau_{j_0}^+$ from Eq.(\ref{tauj}).
\end{description}
Then we have that when $k=k_0$, $\tau=\tau_0$, system (\ref{vdp})
undergoes a Hopf-Hopf bifurcation.

By the above method, we can obtain the Hopf-Hopf bifurcation value,
but estimating the ratio of $\omega_\pm$ is necessary to determine
whether this point is a nonresonant Hopf-Hopf point. Another
algorithm to detect a $k_1:k_2$ resonant Hopf-Hopf point can be
found in \cite{Xujian}. Here we can't use their approach because of
the complexity of the characteristic equation.
\\
\textbf{\textit{Step 2.}}
 Now we will use the algorithm in section 2
to calculate the normal form of (\ref{vdp}) when a Hopf-Hopf
bifurcation occurs at $(k,\tau)=(k_0, \tau_0)$.

When $\tau> 0$, re-scale $t\rightarrow t/\tau$, and denote by
$(k,\tau)=(1/\varepsilon+\alpha_1, \tau_0+\alpha_2)$ we have an
equivalent form of $(\ref{vdp})$:
\begin{equation}\label{vdptaualpha}\left\{\begin{array}{l}\dot x={\left(\tau_0+\alpha_2\right)} y\\
\dot y-\mu\dot
y(t-1)=\left(\tau_0+\alpha_2\right)\{\left[-1+\varepsilon
\left(k_0+\alpha_1\right)(1-\mu)\right]x+\varepsilon y+\mu
x(t-1)\\~~~~~~~~~~~~~-\varepsilon\mu y(t-1)-\varepsilon
x^2y+\varepsilon\mu x^2(t-1)y(t-1)\}\end{array}\right.\end{equation}
When $\alpha_1=\alpha_2=0$, the corresponding characteristic
equation has four roots with zero real parts $\pm i \omega_1\tau_0,
\pm i \omega_2\tau_0$.  Following the procedure in Section 2, we
choose
$$B=\left(
                                                 \begin{array}{cccc}
                                                   i\omega_1\tau_0 & 0 & 0 & 0 \\
                                                   0 & -i\omega_1\tau_0 &  0 & 0 \\
                                                   0 & 0 & i\omega_2\tau_0 & 0 \\
                                                   0 & 0 & 0 & -i\omega_2\tau_0 \
                                                 \end{array}
                                               \right)
,$$

$$\eta(\theta,\alpha)=\left\{
                        \begin{array}{ll}
                          0, & \theta=0; \\
                          -B_1, & \theta\in(-1,0); \\
                          -B_1-B_2, & \theta=-1.
                        \end{array}
                      \right.,
$$
with

$$B_1=\left(
        \begin{array}{cc}
          0 & \tau_0+\alpha_2 \\
          (\tau_0+\alpha_2)[-1+\varepsilon(k_0+\alpha_1)(1-\mu)] & \varepsilon(\tau_0+\alpha_2) \\
        \end{array}
      \right)
$$
and
$$B_2=\left(
        \begin{array}{cc}
          0 & 0 \\(\tau_0+\alpha_2)\mu & -(\tau_0+\alpha_2)\varepsilon\mu \\
        \end{array}
      \right).
$$
Thus we have
$$\Phi(\theta)=\left(
\begin{array}{cccc}
 e^{i \theta  \tau_0  \omega_1 } & e^{-i \theta  \tau_0  \omega_1 } & e^{i
   \theta  \omega_2  \tau_0 } & e^{-i \theta  \omega_2  \tau_0 } \\
 i e^{i \theta  \tau_0  \omega_1 } \omega_1  & -i e^{-i \theta  \tau_0  \omega_1
   } \omega_1  & i e^{i \theta  \omega_2  \tau_0 } \omega_2  & -i e^{-i \theta
   \omega_2  \tau_0 } \omega_2
\end{array}
\right),$$

$$\Psi(s)=\left(
\begin{array}{cc}
 D_1 e^{-i s \tau_0  \omega_1 } \left(-e^{-i \tau_0  \omega_1 } \mu
   \epsilon +\epsilon +i e^{-i \tau_0  \omega_1 } \mu  \omega_1 -i \omega_1
   \right) & -D_1 e^{-i s \tau_0  \omega_1 } \\
 \bar D_1 e^{i s \tau_0  \omega_1 } \left(-e^{i \tau_0  \omega_1 } \mu
   \epsilon +\epsilon -i e^{i \tau_0  \omega_1 } \mu  \omega_1 +i \omega_1
   \right) & -\bar D_1 e^{i s \tau_0  \omega_1 } \\
 D_2 e^{-i s \omega_2  \tau_0 } \left(-e^{-i \omega_2  \tau_0 } \mu
   \epsilon +\epsilon -i \omega_2 +i e^{-i \omega_2  \tau_0 } \omega_2  \mu
   \right) & -D_2 e^{-i s \omega_2  \tau_0 } \\
 \bar D_2 e^{i s \omega_2  \tau_0 } \left(-e^{i \omega_2  \tau_0 } \mu
   \epsilon +\epsilon +i \omega_2 -i e^{i\omega_2 \tau_0 } \omega_2  \mu \right) &
   -\bar D_2 e^{i s \omega_2  \tau_0 }
\end{array}
\right)$$ where $$D_1=\frac{1}{e^{-i \tau_0  \omega_1 } \left(i \mu
\omega_1  (\varepsilon  \tau_0
   +2)-\mu  (\varepsilon +\tau_0 )+\mu  \tau_0  \omega_1 ^2\right)+\varepsilon -2
   i \omega_1 }$$ and $$D_2=\frac{1}{e^{-i \tau_0  \omega_2 } \left(i \mu
\omega_2  (\varepsilon  \tau_0
   +2)-\mu  (\varepsilon +\tau_0 )+\mu  \tau_0  \omega_2 ^2\right)+\varepsilon -2
   i \omega_2 }.$$
\textbf{\textit{Step 3.}} Decomposing  Eq.(\ref{vdptaualpha}) as
Eq.(\ref{decomsimple}),
 we have the form
\begin{equation}\label{deco}\begin{array}{l}
    \dot z_1=i\omega_1\tau_0z_1+D_1 \alpha _2 \left(z_1 \omega _1-z_2 \omega
   _1+\left(z_3-z_4\right) \omega _2\right) \left(\omega _1+i \epsilon
   \right) e^{-i \tau_0  \omega _1} \left(-\mu +e^{i \tau_0  \omega
   _1}\right)\\~~~~~-D_1 [\left(z_1+z_2+z_3+z_4\right) \left(\alpha
   _2+\tau_0 \right) \left(-\epsilon  (\mu -1) \left(k+\alpha
   _1\right)-1\right)\\~~~~~+\left(z_1+z_2+z_3+z_4\right) \tau_0  (k \epsilon
   (\mu -1)+1)\\~~~~~+i \epsilon  \mu  \left(\alpha _2+\tau_0 \right) \left(z_1
   e^{-i \tau_0  \omega _1}+z_2 e^{i \tau_0  \omega _1}+z_3 e^{-i \tau_0
   \omega _2}+z_4 e^{i \tau_0  \omega _2}\right)^2\\~~~~~ \left(z_1 \omega _1
   e^{-i \tau_0  \omega _1}-z_2 \omega _1 e^{i \tau_0  \omega _1}+\omega _2
   e^{-i \tau_0  \omega _2} \left(z_3-z_4 e^{2 i \tau_0  \omega
   _2}\right)\right)\\~~~~~-i \epsilon  \mu \alpha _2
   \left(z_1 \omega _1 e^{-i \tau_0  \omega _1}-z_2 \omega _1 e^{i \tau_0
   \omega _1}+\omega _2 e^{-i \tau_0  \omega _2} \left(z_3-z_4 e^{2 i \tau_0
    \omega _2}\right)\right)\\~~~~~-i \left(z_1+z_2+z_3+z_4\right)^2 \epsilon
    \left(z_1 \omega _1-z_2 \omega _1+\left(z_3-z_4\right) \omega
   _2\right) \left(\alpha _2+\tau_0 \right)\\~~~~~+i \epsilon  \left(z_1 \omega
   _1-z_2 \omega _1+\left(z_3-z_4\right) \omega _2\right) \alpha
   _2\\~~~~~+\mu  \alpha _2 \left(z_1 e^{-i
   \tau_0  \omega _1}+z_2 e^{i \tau_0  \omega _1}+z_3 e^{-i \tau_0  \omega
   _2}+z_4 e^{i \tau_0  \omega _2}\right)],\\
\dot z_2=-i\omega_1\tau_0z_2+\bar D_1 \alpha _2 \left(z_1 \omega
_1-z_2 \omega
   _1+\left(z_3-z_4\right) \omega _2\right) \left(\omega _1-i \epsilon
   \right) \left(-1+\mu  e^{i \tau_0  \omega _1}\right)\\~~~~~-\bar D_1
   [\left(z_1+z_2+z_3+z_4\right) \left(\alpha _2+\tau_0 \right)
   \left(-\epsilon  (\mu -1) \left(k+\alpha
   _1\right)-1\right)\\~~~~~+\left(z_1+z_2+z_3+z_4\right) \tau_0  (k \epsilon
   (\mu -1)+1)\\~~~~~+i \epsilon  \mu  \left(\alpha _2+\tau_0 \right) \left(z_1
   e^{-i \tau_0  \omega _1}+z_2 e^{i \tau_0  \omega _1}+z_3 e^{-i \tau_0
   \omega _2}+z_4 e^{i \tau_0  \omega _2}\right)^2\\~~~~~ \left(z_1 \omega _1
   e^{-i \tau_0  \omega _1}-z_2 \omega _1 e^{i \tau_0  \omega _1}+\omega
   _2 e^{-i \tau_0  \omega _2} \left(z_3-z_4 e^{2 i \tau_0  \omega
   _2}\right)\right)\\~~~~~-i \epsilon  \mu   \alpha _2
   \left(z_1 \omega _1 e^{-i \tau_0  \omega _1}-z_2 \omega _1 e^{i \tau_0
   \omega _1}+\omega _2 e^{-i \tau_0  \omega _2} \left(z_3-z_4 e^{2 i
   \tau_0  \omega _2}\right)\right)\\~~~~~-i \left(z_1+z_2+z_3+z_4\right){}^2
   \epsilon  \left(z_1 \omega _1-z_2 \omega _1+\left(z_3-z_4\right)
   \omega _2\right) \left(\alpha _2+\tau_0 \right)\\~~~~~+i \epsilon  \left(z_1
   \omega _1-z_2 \omega _1+\left(z_3-z_4\right) \omega _2\right)
   \left(\alpha _2+\tau_0 \right)\\~~~~~+\mu  \alpha _2
   \left(z_1 e^{-i \tau_0  \omega _1}+z_2 e^{i \tau_0  \omega _1}+z_3
   e^{-i \tau_0  \omega _2}+z_4 e^{i \tau_0  \omega _2}\right)],\\
   \dot z_3=i\omega_2\tau_0 z_3+D_2 \alpha _2 \left(z_1 \omega _1-z_2 \omega
   _1+\left(z_3-z_4\right) \omega _2\right) \left(\omega _2+i \epsilon
   \right) e^{-i \tau_0  \omega _2} \left(-\mu +e^{i \tau_0  \omega
   _2}\right)\\~~~~~-D_2 [\left(z_1+z_2+z_3+z_4\right) \left(\alpha
   _2+\tau_0 \right) \left(-\epsilon  (\mu -1) \left(k+\alpha
   _1\right)-1\right)\\~~~~~+\left(z_1+z_2+z_3+z_4\right) \tau_0  (k \epsilon
   (\mu -1)+1)\\~~~~~+i \epsilon  \mu  \left(\alpha _2+\tau_0 \right) \left(z_1
   e^{-i \tau_0  \omega _1}+z_2 e^{i \tau_0  \omega _1}+z_3 e^{-i \tau_0
   \omega _2}+z_4 e^{i \tau_0  \omega _2}\right)^2 \\~~~~~ \left(z_1 \omega _1
   e^{-i \tau_0  \omega _1}-z_2 \omega _1 e^{i \tau_0  \omega _1}+\omega _2
   e^{-i \tau_0  \omega _2} \left(z_3-z_4 e^{2 i \tau_0  \omega
   _2}\right)\right)\\~~~~~-i \epsilon  \mu  \alpha _2
   \left(z_1 \omega _1 e^{-i \tau_0  \omega _1}-z_2 \omega _1 e^{i \tau_0
   \omega _1}+\omega _2 e^{-i \tau_0  \omega _2} \left(z_3-z_4 e^{2 i \tau_0
    \omega _2}\right)\right)\\~~~~~-i \left(z_1+z_2+z_3+z_4\right)^2 \epsilon
    \left(z_1 \omega _1-z_2 \omega _1+\left(z_3-z_4\right) \omega
   _2\right) \left(\alpha _2+\tau_0 \right)\\~~~~~+i \epsilon  \left(z_1 \omega
   _1-z_2 \omega _1+\left(z_3-z_4\right) \omega _2\right)  \alpha
   _2 \\~~~~~+\mu  \alpha _2  \left(z_1 e^{-i
   \tau_0  \omega _1}+z_2 e^{i \tau_0  \omega _1}+z_3 e^{-i \tau_0  \omega
   _2}+z_4 e^{i \tau_0  \omega _2}\right)]\end{array}
\end{equation}
$$\begin{array}{l} \dot z_4=-i\omega_2\tau_0z_4+\bar D_2 \alpha _2 \left(z_1 \omega _1-z_2 \omega
   _1+\left(z_3-z_4\right) \omega _2\right) \left(\omega _2-i \epsilon
   \right) \left(-1+\mu  e^{i \tau_0  \omega _2}\right)\\~~~~~-\bar D_2
  [\left(z_1+z_2+z_3+z_4\right) \left(\alpha _2+\tau_0 \right)
   \left(-\epsilon  (\mu -1) \left(k+\alpha
   _1\right)-1\right)\\~~~~~+\left(z_1+z_2+z_3+z_4\right) \tau_0  (k \epsilon
   (\mu -1)+1)\\~~~~~+i \epsilon  \mu  \left(\alpha _2+\tau_0 \right) \left(z_1
   e^{-i \tau_0  \omega _1}+z_2 e^{i \tau_0  \omega _1}+z_3 e^{-i \tau_0
   \omega _2}+z_4 e^{i \tau_0  \omega _2}\right)^2\\~~~~~ \left(z_1 \omega _1
   e^{-i \tau_0  \omega _1}-z_2 \omega _1 e^{i \tau_0  \omega _1}+\omega
   _2 e^{-i \tau_0  \omega _2} \left(z_3-z_4 e^{2 i \tau_0  \omega
   _2}\right)\right)\\~~~~~-i \epsilon  \mu   \alpha _2
   \left(z_1 \omega _1 e^{-i \tau_0  \omega _1}-z_2 \omega _1 e^{i \tau_0
   \omega _1}+\omega _2 e^{-i \tau_0  \omega _2} \left(z_3-z_4 e^{2 i
   \tau_0  \omega _2}\right)\right)\\~~~~~-i \left(z_1+z_2+z_3+z_4\right)^2
   \epsilon  \left(z_1 \omega _1-z_2 \omega _1+\left(z_3-z_4\right)
   \omega _2\right) \left(\alpha _2+\tau_0 \right)\\~~~~~+i \epsilon  \left(z_1
   \omega _1-z_2 \omega _1+\left(z_3-z_4\right) \omega _2\right)
   \alpha _2 \\~~~~~+\mu   \alpha _2
   \left(z_1 e^{-i \tau_0  \omega _1}+z_2 e^{i \tau_0  \omega _1}+z_3
   e^{-i \tau_0  \omega _2}+z_4 e^{i \tau_0  \omega _2}\right)].\end{array}$$
   Following the algorithm in section 2 and  doing the projection of Eq.(\ref{deco}) onto $(Im(M_2^1))^c$ and $(Im(M_3^1))^c$ then we
   have these coefficients
\begin{eqnarray}\label{ac}
  \nonumber a_{11} &=& -D_1 \varepsilon  (1-\mu ) \tau_0 \\
 \nonumber  a_{12} &=& D_1 \left(k_0 \varepsilon  (\mu -1)-\mu  \left(\omega_1 ^2+1\right)
   e^{-i \tau_0  \omega_1 }+\omega_1 ^2+1\right)\\
 \nonumber  c_{11} &=& -\frac{1}{2} D_1 \left(2 i \varepsilon  \mu  \tau_0  \omega_1  e^{-i
   \tau_0  \omega_1 }-2 i \varepsilon  \tau_0  \omega_1 \right) \\
\nonumber  c_{12} &=& -D_1 \left(2 i \varepsilon  \mu  \tau_0
\omega_1 e^{-i \tau_0  \omega_1
   }-2 i \varepsilon  \tau_0  \omega_1 \right) \\
   \nonumber a_{21} &=& -D_2 \varepsilon  (1-\mu ) \tau_0 \\
 \nonumber  a_{22} &=& D_2 \left(k_0 \varepsilon  (\mu -1)-\mu  \left(\omega_2 ^2+1\right)
   e^{-i \tau_0  \omega_2 }+\omega_2 ^2+1\right)\\
 \nonumber  c_{21} &=&  -D_2 \left(2 i \varepsilon  \mu  \tau_0 \omega_2
e^{-i \tau_0  \omega_2
   }-2 i \varepsilon  \tau_0  \omega_2 \right)\\
\nonumber  c_{22} &=&-\frac{1}{2} D_2 \left(2 i \varepsilon  \mu
\tau_0  \omega_2  e^{-i
   \tau_0  \omega_2 }-2 i \varepsilon  \tau_0  \omega_2 \right)  \\
\end{eqnarray}
Substituting Eq.(\ref{ac}) into Eq.(\ref{array}), we can distinguish
the unfoldings by Table 1. So far, all the key coefficients
determining the normal form in Eq. (\ref{NF23polar}) are obtained.
However, due to the complexity of the van der Pol's equation, it is
quite difficult to estimate the sign of $b_0,~c_0,~d_0$, and
$d_0-b_0c_0$, thus we give a numerical example in the coming
section.
\subsection{Illustrations}
\begin{figure}
  \begin{center}\includegraphics[width=8cm]{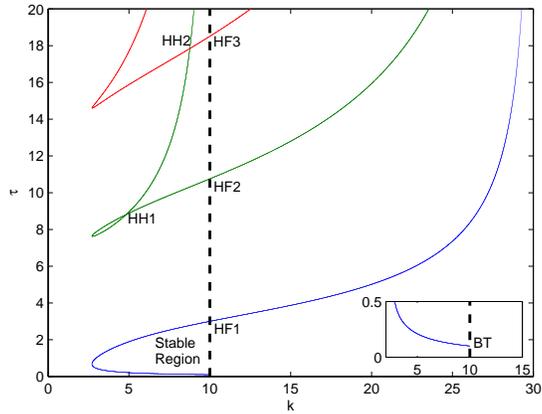}\end{center}
  \caption{Partial bifurcation sets with parameters in the $k-\tau$ plane. The color lines are Hopf bifurcation curve and the dotted line stands for the fold bifurcation curve.}\label{fig1}
\end{figure}

\begin{figure}
  \begin{center}\includegraphics[width=8cm]{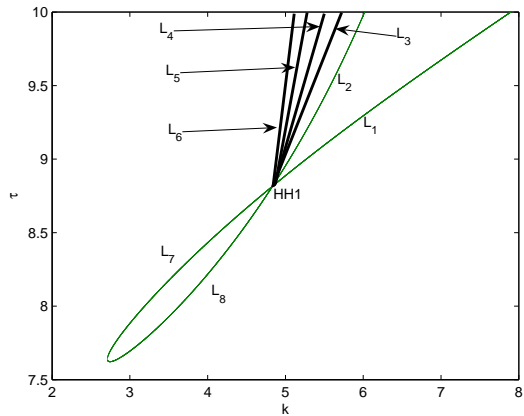}\end{center}
  \caption{Complete bifurcation sets near HH1.}\label{fig2}
\end{figure}

\begin{figure}
  \begin{center}\includegraphics[width=13cm]{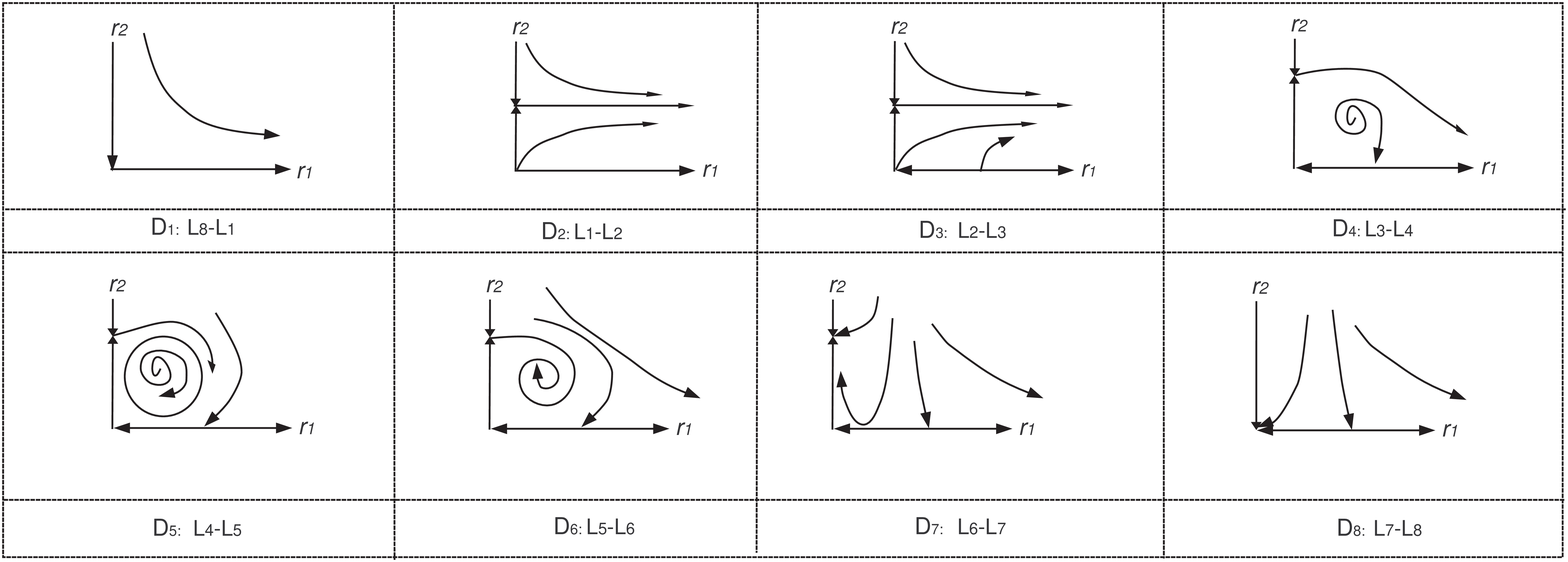}\end{center}
  \caption{The eight distinct phase portraits near HH1 in $D_1$--$D_8$. Below every figure, we mark the corresponding region in Figure \ref{fig2}, e.g. the region $D_1$ is between $L_8$ and $L_1$.}\label{fig3}
\end{figure}

In this section we choose $\varepsilon=0.1$.  Adding the extended
delay feedback with $\mu=0.5$ into (\ref{vdp}), following the
regular characteristic equation analysis and Theorem 2, we have the
bifurcation diagram in the $k-\tau$ plane as in Figure \ref{fig1}.
Here we only state the main results about Hopf bifurcation. In
Figure \ref{fig1}, several colored Hopf bifurcation curves and a
dotted fold bifurcation curve are presented. When $\tau=0$ the zero
solution is unstable and a stable region of the zero solution is
marked by \lq\lq Stable Region". One Bogdanov-Takens point, three
Hopf-fold points and two Hopf-Hopf points are marked by BT, HF1-HF3
and HH1-HH2, respectively.  From Eq. (\ref{cossin}), (\ref{rho}) and
(\ref{k0}) we have when $$k_0 = 4.834585253, ~\tau_0 =
8.815987316,$$ two different frequencies are solved  by $$\omega_1 =
0.7307969965,$$ and
$$\omega_2 = 0.90073546761$$ with
$\omega_1:\omega_2=0.811334:1$, thus this point is a nonresonant
Hopf-Hopf bifurcation point. Following Eq.(\ref{array}) and
Eq.(\ref{ac}) we have
$$c_1=0.2429777596\alpha_1-0.2981855434\alpha_2,$$
$$c_2=-0.2004123093\alpha_1+0.4602126544\alpha_2,$$
$$b_0=0.087454,$$ $$c_0=-45.7383,$$ $$d_0=-1,$$ and $$d_0-b_0c_0=3.$$
By Table 1, we know the case VIa arises. From Guckenheimer
\cite{Guckenheimer}, near the Hopf-Hopf point $HH1$ there are eight
different  kinds of phase diagrams in eight different regions which
are divided by lines
 $L_1$--$L_8$ with
\begin{description}
\item [$L_1$]: $\alpha_2=0.435478\alpha_1,~\alpha_1>0$;
\item [$L_2$]: $\alpha_2=0.814854\alpha_1,~\alpha_1>0$;
\item [$L_3$]: $\alpha_2=0.828102\alpha_1,~\alpha_1>0$;
\item [$L_4$]: $\alpha_2=0.828985\alpha_1+O(\alpha_1^2),~\alpha_1>0$;
\item [$L_5$]: $\alpha_2=0.828985\alpha_1,~\alpha_1>0$;
\item [$L_6$]: $\alpha_2=0.874050\alpha_1,~\alpha_1>0$;
\item [$L_7$]: $\alpha_2=0.435478\alpha_1,~\alpha_1<0$;
\item [$L_8$]: $\alpha_2=0.814854\alpha_1,~\alpha_1<0$;
\end{description}
Recall that $$\alpha_1=k-k_0,~\alpha_2=\tau-\tau_0,$$ thus we give a
bifurcation set on the plane of the original parameters in system
(\ref{vdp}) (See Figure \ref{fig2}). In figure \ref{fig3}, we draw
these phase portraits and label the position where the corresponding
parameters lie in. In every portrait, a nontrivial equilibrium on
the axis,  an equilibrium with positive $r_1,~r_2$  and a cycle
correspond to a nonconstant periodic, a quasi-periodic solution on
the 2-dimension torus and a quasi-periodic  solution on the
3-dimension torus of Eq.(\ref{vdp}), respectively.

\begin{figure}
  \begin{center}\includegraphics[width=8cm]{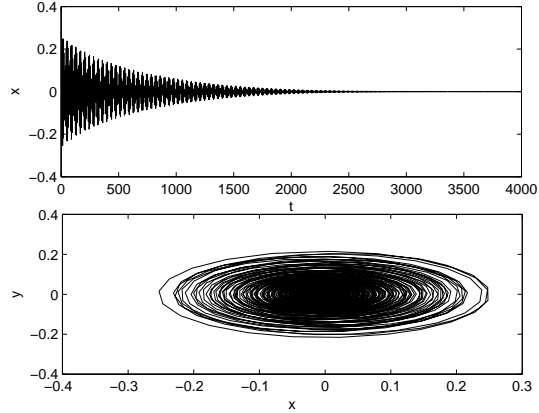}\end{center}
  \caption{$\alpha_1=-0.1$, $\alpha_2=-0.08$ in $D_8$. The trivial equilibrium of system (\ref{vdp}) is stable.}\label{fig4}
\end{figure}

Now we give some simulations. When $\alpha_1=-0.1,~ \alpha_2=-0.08$
(in Region $D_8$), system (\ref{vdp}) has a stable equilibrium,
which is shown in Figure \ref{fig4}.

\begin{figure}
  \begin{center}\includegraphics[width=9cm]{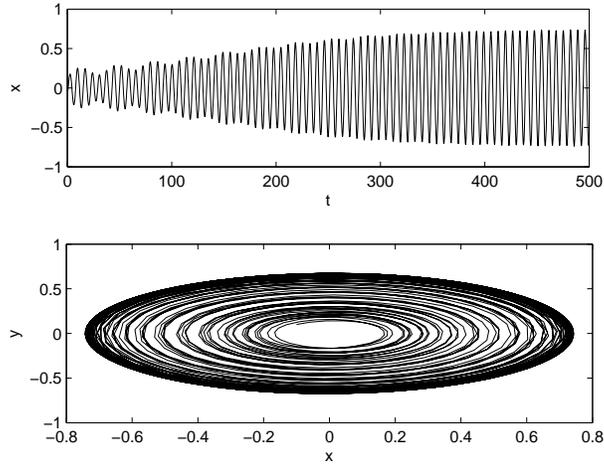}\end{center}
  \caption{$\alpha_1=-0.1$, $\alpha_2=0.1$ in $D_7$. System (\ref{vdp}) has a stable periodic solution in Region $D_7$.}\label{fig5}
\end{figure}

In $D_7$, Figure \ref{fig3} indicates there is a stable periodic
solution, which is also illustrated in Figure \ref{fig5}, where
$\alpha_1=-0.1$, $\alpha_2=0.1$.

\begin{figure}
  \begin{center}\includegraphics[width=9cm]{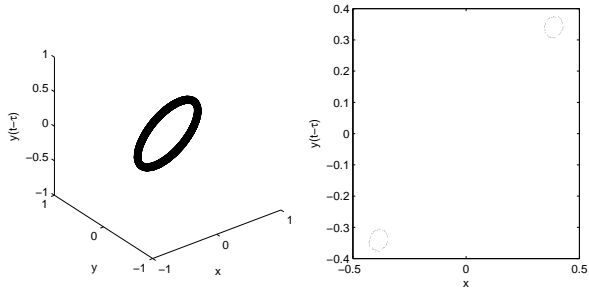}\end{center}
  \caption{$\alpha_1=0.1$, $\alpha_2=0.085$ in $D_6$. The bifurcated quasi-periodic solution of system (\ref{vdp}). }\label{fig6}
\end{figure}

When parameters are chosen between $L_5$ and $L_6$ (i.e. in $D_6$),
there exists a stable quasi-periodic solution on a 2-dimensional
torus which is shown in Figure \ref{fig6}, where $\alpha_1=0.1$,
$\alpha_2=0.085$.

When parameters are chosen between $L_4$ and $L_5$ (i.e. in $D_5$),
there exists a quasi-periodic solution on a 3-dimensional torus
which is shown in Figure \ref{fig6}, where $\alpha_1=0.2$,
$\alpha_2=0.164$. The right figure is the
Poincar$\acute{\textrm{e}}$ map on the whole
Poincar$\acute{\textrm{e}}$ section $y(t)=0$. Clearly, we find that
the points on the Poincar$\acute{\textrm{e}}$ section exhibits
quasi-periodic behavior, which indicates the solution is a
quasi-periodic solution on a 3-dimensional torus.

\begin{figure}
  \begin{center}\includegraphics[width=13cm]{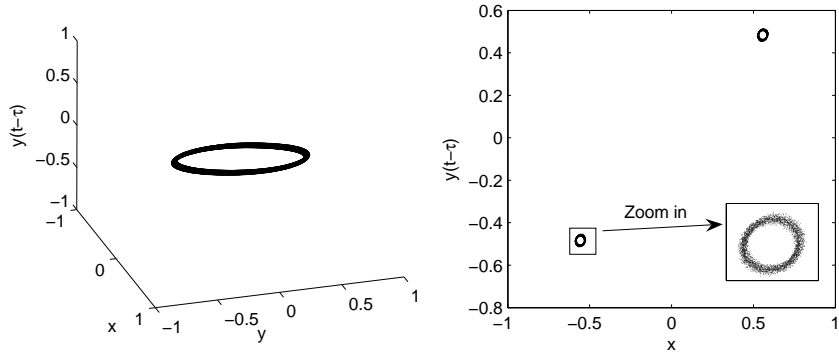}\end{center}
  \caption{$\alpha_1=0.2$, $\alpha_2=0.164$ in $D_5$. The bifurcated quasi-periodic solution on the three-dimensional torus of system (\ref{vdp}) and the corresponding Poincar$\acute{\textrm{e}}$ map on the whole Poincar$\acute{\textrm{e}}$ section $y(t)=0$. }\label{fig7}
\end{figure}

\begin{figure}
  \begin{center}\includegraphics[width=7cm]{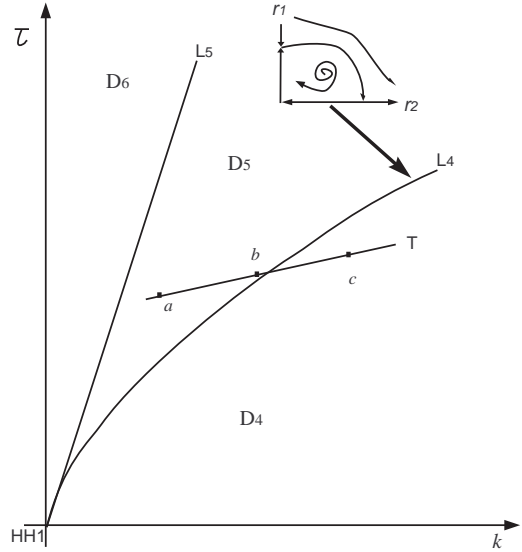}\end{center}
  \caption{ The sketch of the saddle connection bifurcation curve ($L_4$) on the three-dimensional torus. (a) $(\alpha_1, \alpha_2)=2\times(0.1, 0.081)$; (b) $(\alpha_1, \alpha_2)=2.4\times(0.1, 0.081)$; (c) $(\alpha_1, \alpha_2)=2.5\times(0.1, 0.081)$.  }\label{fig8}
\end{figure}

\begin{figure}
  \begin{center}\includegraphics[width=14cm]{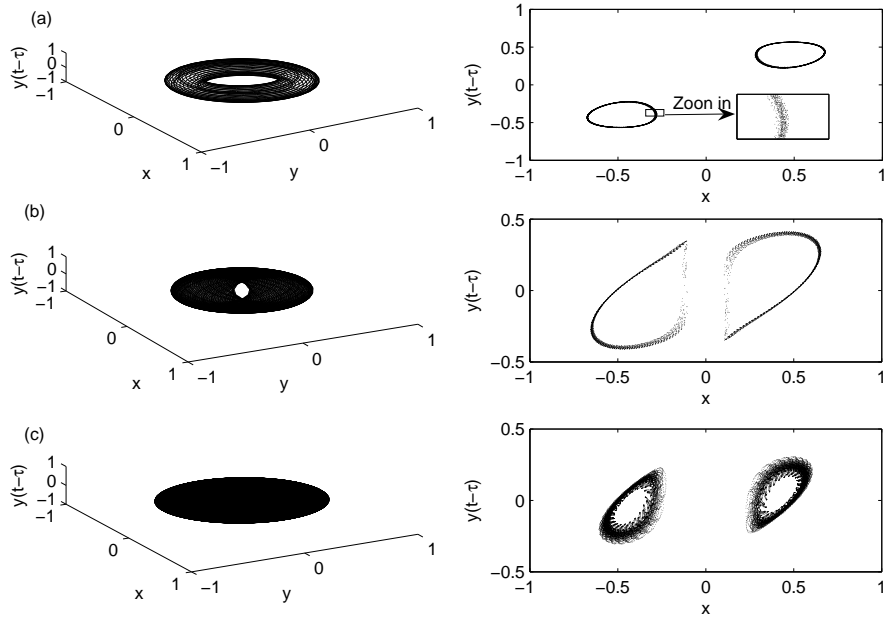}\end{center}
  \caption{$\alpha_1=0.2$, $\alpha_2=0.164$. The phase portraits in $x$--$y$--$y(\cdot-\tau)$ space and the corresponding Poincar$\acute{\textrm{e}}$ map on the whole Poincar$\acute{\textrm{e}}$ section $y(t)=0$ when parameters are chosen at (a), (b) and (c) in Figure \ref{fig8}, respectively. }\label{fig9}
\end{figure}

\begin{figure}
  \begin{center}\includegraphics[width=14cm]{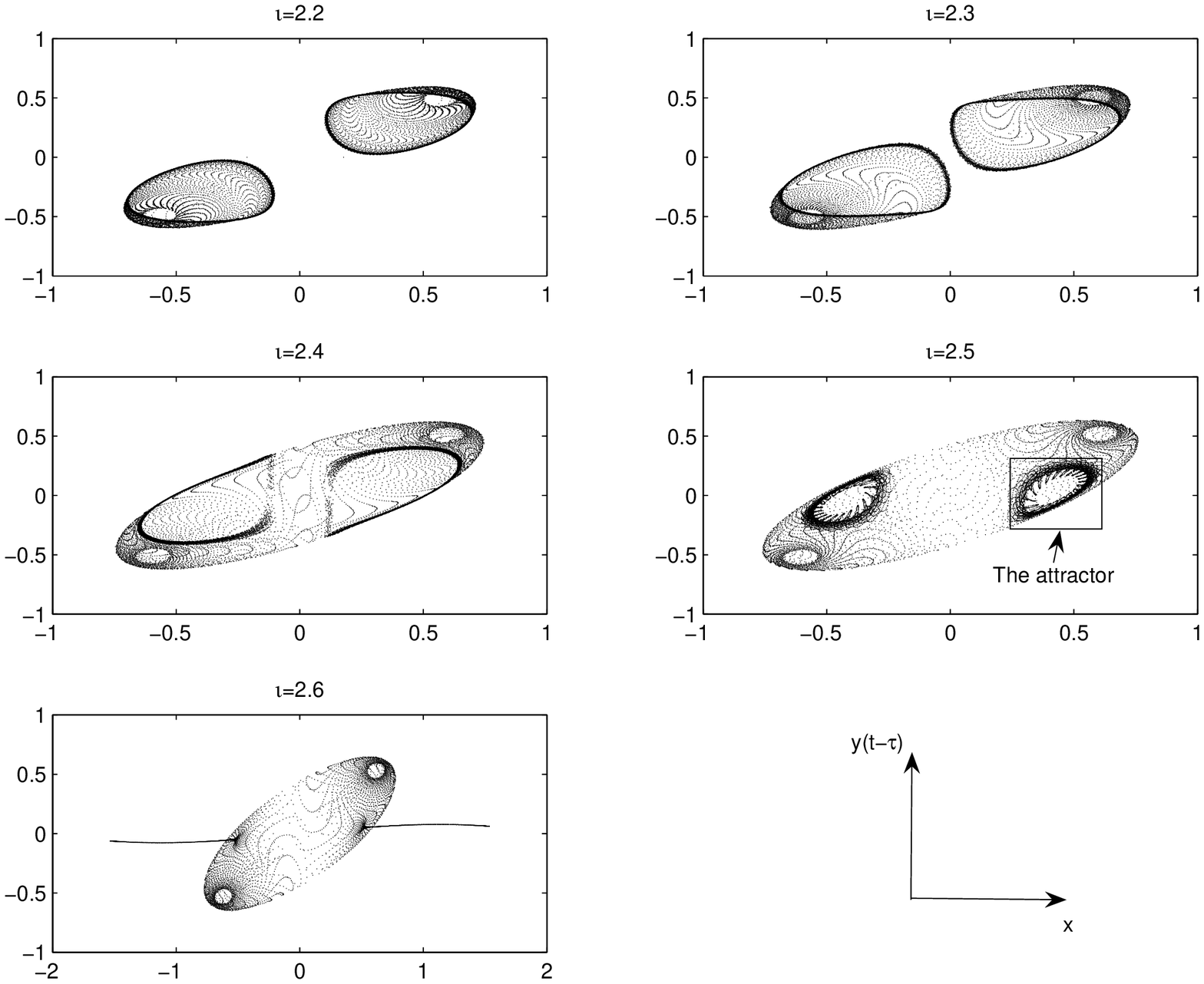}\end{center}
  \caption{The complete Poincar\'e maps near the saddle connection.}\label{fig10}
\end{figure}

Generally, a vanishing 3-dimensional torus might bring chaos to the
system\cite{Battelino,eck,rue}. In system (\ref{vdp}), we choose three points (a), (b) and
(c) on the line T: $(\alpha_1,\alpha_2)=(0.1\iota, 0.081\iota)$,
which is shown in Figure \ref{fig8}. In Figure \ref{fig9}, the phase
portraits are drawn. At (a), the system has a quasi-periodic
solution on a three-dimensional torus, which vanishes via the saddle
connection bifurcation on the three-dimensional torus (the curve
$L_4$). At point (c), system (\ref{vdp}) exhibits chaotic behavior
and the chaotic attractor is drawn on the
Poincar$\acute{\textrm{e}}$ section $y(t)=0$. Thus we confirm that
in NDDE the destroying of a three-dimensional torus might bring
chaos as mentioned in \cite{Battelino}. In order to give a neat
expression we delete the transient states in Figure \ref{fig7} and
\ref{fig9}. Figure \ref{fig10} is also an illustration of the transition where we give the complete Poincar\'e map, from which we find the strange attractor vanishing when $\iota=2.6$. After that, the system is stabilized to a
periodic solution with large amplitude.

\section{Conclusions}
In this paper we mainly study the nonresonant Hopf-Hopf bifurcation
in a NFDE with parameters as follows
\begin{equation}\label{NFDEc}
    \frac d {dt}\left[Dx_t-G(x_t)\right]=L(\alpha)x_t+F(\alpha, x_t)
\end{equation}
Following \cite{Faria1, Faria2, Wee, WangWei1}, we compute the
normal form near the bifurcation point. An explicit algorithm is
given  to calculate the four key variables: $b_0,~c_0,~d_0$ and
$d_0-b_0c_0$, by which the twelve unfoldings are distinguished. We
find that the operator $D$ just changes the method when transforming
the NFDE to a abstract ODE and the decomposing of phase space,
compared with the normal form derivation for RFDE. All the rest
steps of calculations remain almost the same as dealing with a RFDE.

As an illustration of this theory, van der Pol's equation with
extended delay feedback is considered. We give the conditions under
which the Hopf-Hopf bifurcation occurs.  Detailed dynamics near the
origin are obtained by drawing the corresponding bifurcation set.
Both theoretical bifurcation set and simulations confirm the
existence of stable periodic solutions and stable quasi-periodic
solutions. With the guide of the bifurcation sets we also find in
van der Pol's equation a chaotic attractor appears as the
three-dimensional torus vanishes via a saddle connection
bifurcation.







\end{document}